\documentclass{amsart}
\usepackage{amssymb,amscd}
\usepackage{graphicx}

\usepackage{latexsym}

\usepackage{amsthm}
\newtheorem{theorem}{Theorem}
\newtheorem{condition}{Condition}
\newtheorem*{defn}{Definition}
\newtheorem*{lemma}{Lemma}
\newtheorem*{cor}{Corollary}
\newtheorem{prop}{Proposition}[section]

\newcommand{\bc}{\mathbb{C}}
\newcommand{\bh}{\mathbb{H}}
\newcommand{\bp}{\mathbb{P}}

\newcommand{\br}{\mathbb{R}}
\newcommand{\bz}{\mathbb{Z}}
\newcommand{\modm}{\mathcal{M}}
\newcommand{\modt}{\mathcal{T}}

\newcommand{\tr}{{\rm tr}\hspace{.2mm}}
\newcommand{\pants}{\includegraphics[height=.4cm]{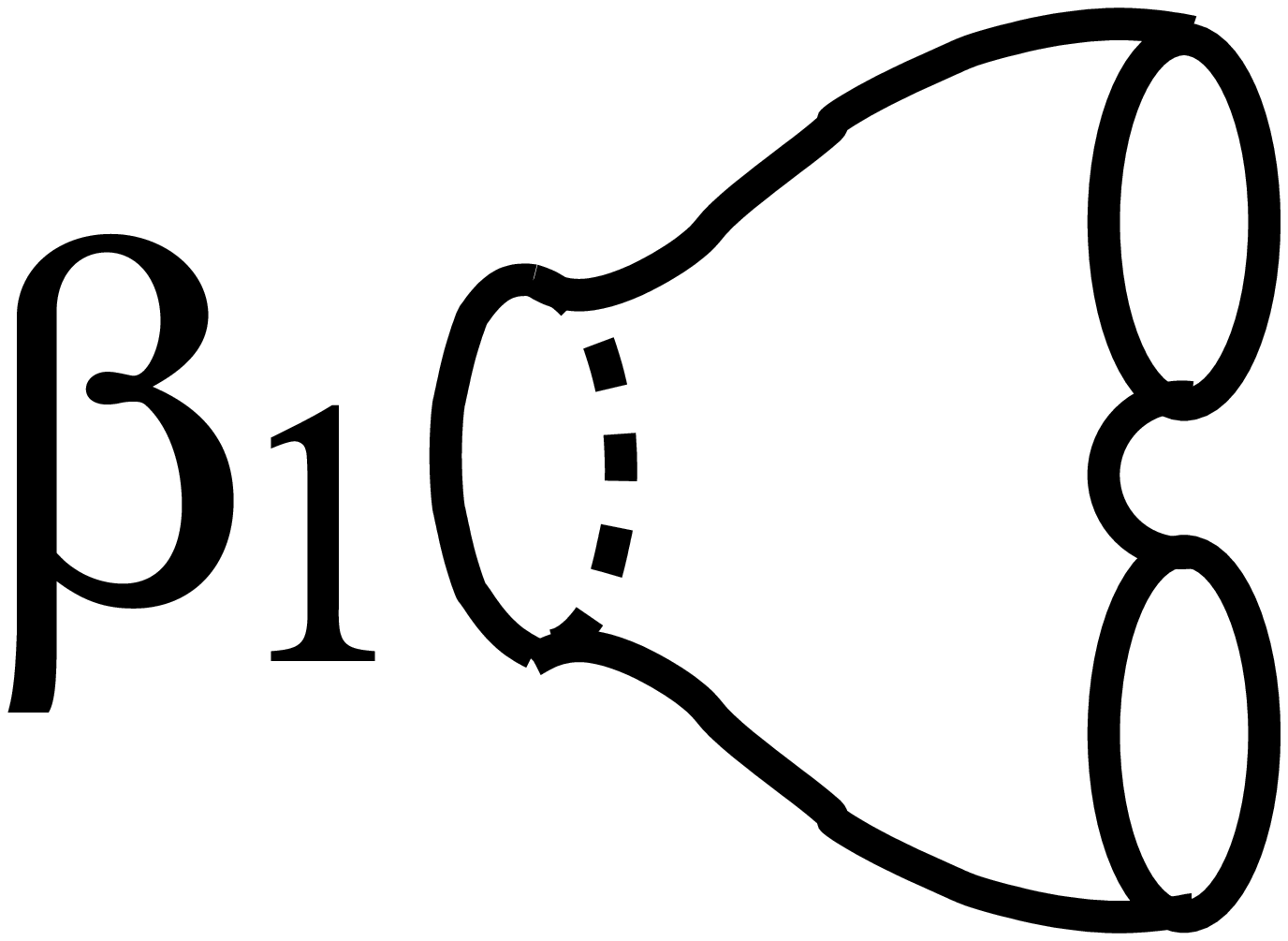}}
\newcommand{\mobius}{\includegraphics[height=.4cm]{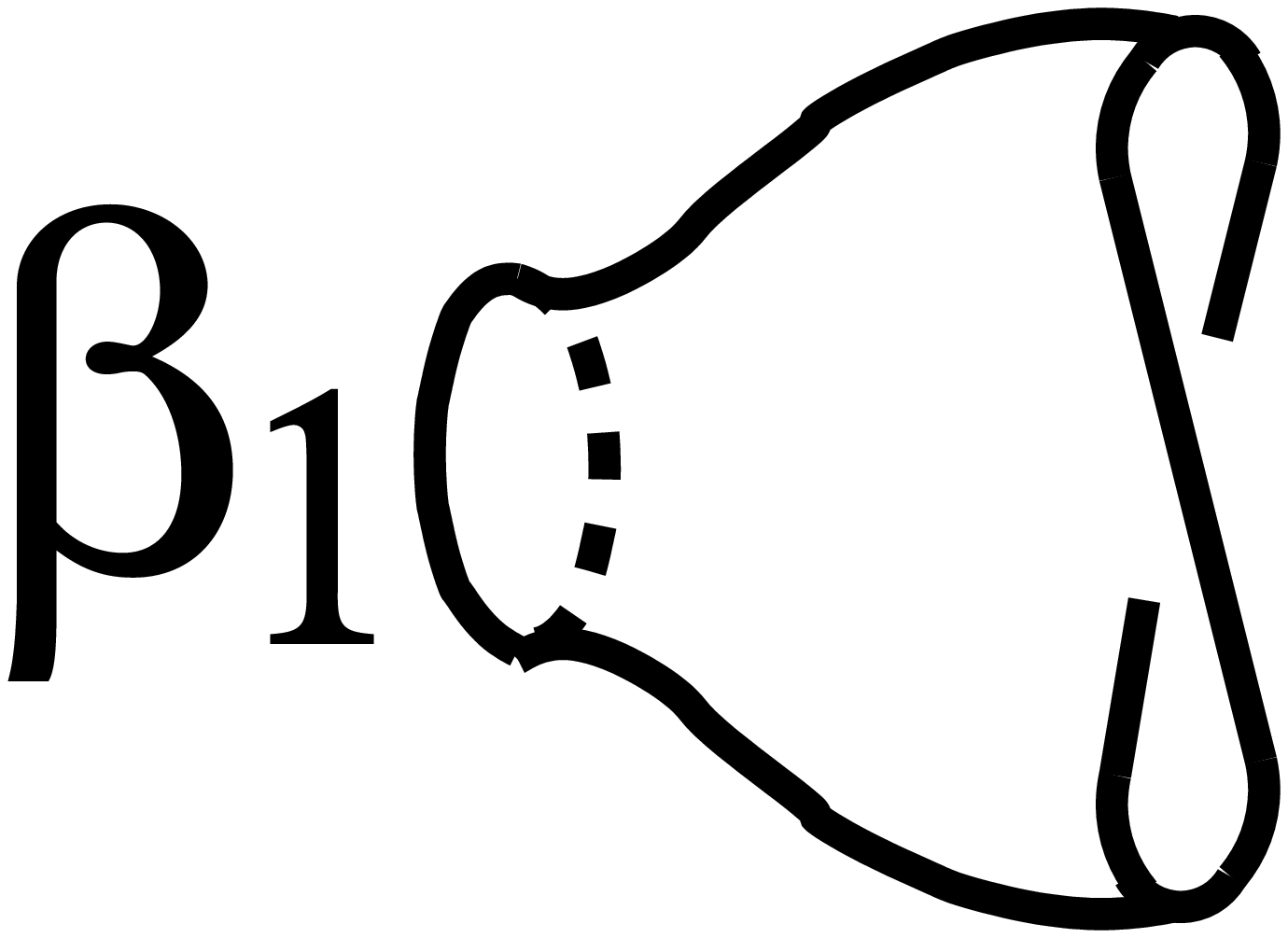}}

\begin{document}

\title{Lengths of geodesics on non-orientable hyperbolic surfaces}
\author{Paul Norbury}
\address{Department of Mathematics and Statistics\\
University of Melbourne\\Australia 3010\\
and Boston University\\
111 Cummington St\\
Boston MA 02215}
\email{pnorbury@math.bu.edu}
\thanks{The author would like to thank the Max Planck Institut F\"ur Mathematik, Bonn where part of this work was carried out.}

\date{\today}

\begin{abstract}

\noindent We give an identity involving sums of functions of lengths of simple closed geodesics, known as a McShane identity, on any non-orientable hyperbolic surface with boundary which generalises Mirzakhani's identities on orientable hyperbolic surfaces with boundary.

\end{abstract}

\maketitle

\section{Introduction}

Closed geodesics on hyperbolic surfaces have extremely rich properties, arising in geometry, topology and number theory.  On a hyperbolic surface the set of closed geodesics of length less than any constant is finite so in particular the set of all closed geodesics is countable.  One can sum functions that are sufficiently decreasing, of lengths of geodesics, over the set of all closed geodesics, and the most famous identity involving such a sum is the Selberg trace formula.   An identity involving sums of functions of lengths of simple (embedded) closed geodesics on a punctured hyperbolic torus was obtained by McShane \cite{McSRem}, whereas the Selberg trace formula includes the non-embedded closed geodesics.  McShane's identity was generalised by various authors \cite{McSSim,MirSim,TWZGen} to other orientable hyperbolic surfaces.  Most notably, in \cite{MirSim} Mirzakhani applied the identities to get deep information about the moduli space of hyperbolic surfaces.

In this paper we prove a McShane identity---a non-trivial sum of functions of lengths of geodesics---on non-orientable hyperbolic surfaces with geodesic boundary. This was motivated from two directions.  Firstly, Mirzakhani \cite{MirWei} used her generalisations of McShane's identity to study recursion relations among tautological cohomology classes on the moduli spaces.  Recently Wahl \cite{WahHom} has shown that the tautological cohomology classes on moduli spaces of non-orientable hyperbolic surfaces have stable behaviour analogous to the orientable case.  One would expect there also to be recursion relations among these cohomology classes and that such recursion relations may interact with the known recursion relations.  The McShane identity proven in this paper is potentially a tool to find such relations.  Secondly, orientable hyperbolic surfaces with geodesic boundary lie in a larger class of surfaces known as Klein surfaces.  A {\em Klein surface} is a real algebraic curve, meaning it is a complex algebraic curve with only real coefficients in its defining equation.  It possesses an anti-holomorphic involution and can be identified with its quotient by the involution.  For example, inside the moduli space of elliptic curves, identified as a fundamental domain of $PSL(2,\bz)$ acting on upper half-space $\bh^2$, the Klein surfaces correspond to curves with representatives $z\in\bh^2$ satisfying $Re(z)\in\frac{1}{2}\bz$.  With respect to the unique hyperbolic structure on a Klein surface, the anti-holomorphic involution is an isometry with fixed point set a (possibly empty) collection of simple closed geodesics.  The quotient is a hyperbolic surface with geodesic boundary that corresponds to the fixed point set of the involution.  It is natural to try to extend the McShane identities on orientable hyperbolic surfaces with geodesic boundary to all Klein surfaces.  We extend the identities to Klein surfaces with involution having non-empty fixed point set, equivalently the quotient has non-empty boundary, leaving only those with involutions that act freely.

We use the following standard terminology to describe the global and local behaviour of a closed curve on a surface.  A {\em simple} closed curve is an embedded curve.  A simple closed curve is either globally {\em separating} or {\em non-separating}, and it is locally  {\em one-sided} or {\em two-sided}.   As usual, we omit the words globally and locally in the paper.  Notice that there are exactly three types of simple closed curves since a one-sided closed curve is necessarily non-separating.

We first state the simplest version of the result, when the hyperbolic surface is the punctured Klein bottle $K$, since it is the most readable form and has independent interest.  Up to isotopy there exists a unique two-sided simple closed curve $\gamma\subset K$  and an infinite sequence of one-sided simple closed curves $\gamma_i$.  The two-sided and one-sided isotopy classes are orbits of the mapping class group of $K$ (which is $\bz\rtimes\bz_2$.)  In other words, the isotopy class of the two-sided simple closed curve is invariant under the mapping class group of $K$ and the isotopy classes of any two one-sided simple closed curves are related by an element of the mapping class group.

Equip $K$ with a complete hyperbolic metric with a cusp.  Denote the length of the unique geodesic in each isotopy class by $l_{\gamma}$ and $l_{\gamma_i}$.  

\begin{theorem}  \label{th:kident}
On any hyperbolic punctured Klein bottle,
\begin{equation}  \label{eq:kident}
\sum_{i=-\infty}^{\infty}\frac{1}{1+\sinh^2\frac{1}{2}l_{\gamma_i}+\sinh^2\frac{1}{2}l_{\gamma_{i+1}}}=\tanh\frac{1}{2}l_{\gamma}
\end{equation}
where $\gamma$ is the unique two-sided geodesic and the sum is over all pairs of disjoint simple closed geodesics.
\end{theorem}
Note that the geodesic $\gamma$, and hence $\tanh(l_{\gamma})$, is well-defined whereas to make sense of the $i$th one-sided simple closed curve $\gamma_i$ requires a marking of $K$.  Nevertheless, the left hand side of (\ref{eq:kident}) is well-defined using the property that $\gamma_i$ is disjoint from both $\gamma_{i-1}$ and $\gamma_{i+1}$ and has non-trivial intersection with all other $\gamma_j$.

\begin{theorem}  \label{th:main}
For
\begin{equation}
R(x,y,z)=x-\ln\frac{\cosh\frac{y}{2}+\cosh\frac{x+z}{2}}{\cosh\frac{y}{2}+\cosh\frac{x-z}{2}}
\end{equation}
and
\[ D(x,y,z)=R(x,y,z)+R(x,z,y)-x,\quad E(x,y,z)=R(x,2z,y)-\frac{x}{2}\]
on a hyperbolic surface with Euler characteristic $\neq-1$ the following identity holds:
\begin{equation}  \label{eq:nonormirz}
\sum_{\gamma_1,\gamma_2}D(L_1,l_{\gamma_1},l_{\gamma_2})+\sum_{j=2}^n\sum_{\gamma}R(L_1,L_j,l_{\gamma})
+\sum_{\mu,\nu}E(L_1,l_{\nu},l_{\mu})=L_1
\end{equation}
where the sums are over simple closed geodesics.   The first sum is over pairs of two-sided geodesics $\gamma_1$ and $\gamma_2$ that bound a pair of pants with $\beta_1$, the second sum is over boundary components $\beta_j$, $j=2,..,n$ and two-sided geodesics $\gamma$ that bound a pair of pants with $\beta_1$ and $\beta_j$, and the third sum is over one-sided geodesics $\mu$ and two-sided geodesics $\nu$ that, with $\beta_1$, bound a M\"{o}bius strip minus a disk containing $\mu$.
\end{theorem}
In the orientable case, the third term in (\ref{eq:nonormirz}) vanishes and we are left with Mirzakhani's identity \cite{MirSim}.  The functions $D$ and $R$ are taken from \cite{MirSim}.  Define $z'(x,y,z)$ by $\cosh\frac{x}{2}+\cosh\frac{y}{2}=2\sinh\frac{z}{2}\sinh\frac{z'}{2}$ (see (\ref{eq:mobident})) then
\[ 0\leq R(x,y,z)\leq x,\quad 0\leq D(x,y,z)\leq x,\quad
0\leq E(x,y,z)+E(x,y,z')\leq x\]
so the quotients of these functions by $x$ are best understood as probabilities, and they make sense even in the limit $x\to 0$.  This is the view taken in Section~\ref{sec:proof}.

The expression $E(x,y,z)+E(x,y,z')$ for $z'(x,y,z)$ defined above arises because the summands $E(L_1,l_{\nu},l_{\mu})$ in (\ref{eq:nonormirz}) naturally come in pairs $E(L_1,l_{\nu},l_{\mu})+E(L_1,l_{\nu},l_{\mu'})$ where $\mu$ and $\mu'$ are the two one-sided geodesics in the M\"{o}bius strip minus a disk bounded by $\nu$ and $\beta_1$ (where $l_{\beta_1}=L_1$).  The property $0\leq E(x,y,z)+E(x,y,z')$ allows us to consider the series as containing only positive summands.  Moreover, we get the necessary decay \[\lim_{y\to\infty}E(x,y,z)+E(x,y,z')=0=\lim_{z\to\infty}E(x,y,z)+E(x,y,z')\]
required for the convergence of the series.

We do not deduce the non-orientable case of Theorem~\ref{th:main} from the orientable case since the proof involves analysing the self-intersections of geodesics on the non-orientable surface, and the orientable double cover does not see these self-intersections in general.

Theorem~\ref{th:main} extends to the four hyperbolic surfaces of Euler characteristic $=-1$, although they each require special treatment.  In Section~\ref{sec:proof} the functions $D(x,y,z)$ and $E(x,y,z)$ are defined independently of $R(x,y,z)$.  The analogue of (\ref{eq:nonormirz}) on a pair of pants produces the relation between $D$ and $R$ given in Theorem~\ref{th:main}, and on the M\"{o}bius strip minus a disk it produces the relation between $E$ and $R$.  

On a Klein bottle minus a disk the McShane identity bares little resemblance to (\ref{eq:nonormirz}) (whereas the relationship on a torus minus a disk in \cite{MirSim} consists of the first term of (\ref{eq:nonormirz}) evaluated on $l_{\gamma_1}=l_{\gamma_2}$.)   Theorem~\ref{th:kident} follows as a limiting case $L\to 0$ of the following theorem.
\begin{theorem}  \label{th:kident2}
On a hyperbolic Klein bottle with geodesic boundary of length $L$, for
\[F(x,y,z)=\frac{x}{2}-\ln\frac{\cosh{y}+\exp\frac{x}{2}\cosh{z}-\sinh\frac{x}{2}}{\cosh{y}+\exp\frac{-x}{2}\cosh{z}+\sinh\frac{x}{2}}\]
\begin{equation}  \label{eq:kident2}
\sum_{\mu,\nu}F(L,l_{\mu},l_{\nu})=2\ln\frac{1+e^{L/2}e^{l_{\gamma}}}{e^{L/2}+e^{l_{\gamma}}}
\end{equation}
where $\gamma$ is the unique two-sided geodesic and the sum is over all ordered pairs of disjoint simple closed geodesics $(\mu,\nu)$.
\end{theorem}

Non-orientable surfaces provide the simplest non-trivial hyperbolic surface---the punctured Klein bottle---which helps our understanding of the general case.  The punctured Klein bottle is simpler than the punctured torus due to the simplicity of its mapping class group---$\bz\rtimes\bz_2$ which can be essentially thought of as its index 2 normal subgroup $\bz$.  In Section~\ref{sec:klein} we give an elementary treatment of this simplest case, which also generalises to complex lengths.  This is analogous to Bowditch's treatment of the punctured torus \cite{BowPro}.

This paper has some orientable consequences.  The representation of a non-orientable hyperbolic surface is an {\em extended quasi-fuchsian} group, or in other words $PGL(2,\br)\subset PSL(2,\bc)$ and non-orientable hyperbolic surfaces give examples of orientable hyperbolic 3-manifolds.  The orientable 3-manifolds are homeomorphic to non-trivial $\br$ bundles over non-orientable surfaces, and admit more general extended quasi-fuchsian groups coming from deforming the hyperbolic surface examples.  The relationship between traces and lengths of geodesics on a hyperbolic surface generalises to a relationship between traces and complex lengths in a hyperbolic 3-manifold.  A {\em complex length} stores two pieces of geometric information.  Its real part is the length of the corresponding closed geodesic, and its imaginary part which is well-defined mod $\pi$ describes the angle of rotation of the normal bundle parallel transported along the closed geodesic.  
\begin{theorem}  \label{th:3man}
On an oriented cusped hyperbolic 3-manifold homeomorphic to a bundle over a punctured Klein bottle
\begin{equation}
\sum_{i=-\infty}^{\infty}\frac{1}{1+\sinh^2(z_i/2)+\sinh^2(z_{i+1}/2)}=\tanh(z/2)
\end{equation}
where $z_i$ and $z_{i+1}$ are the complex lengths of closed geodesics that project onto a pair of disjoint embedded closed curves on the surface and $z$ is the complex length of the unique closed geodesic that projects to an embedded two-sided curve on the surface.
\end{theorem}
\vspace{.5cm}
In the final section, we describe the moduli space of hyperbolic surfaces.  We find that it is natural to consider {\em pairs} consisting of a hyperbolic surface together with an orientation of one of its boundary components.  This is analogous to choosing an orientation on an orientable surface.  It allows the construction of a well-defined volume form on the moduli space which is required for integration over the moduli space.

We finish by describing the consequences that (\ref{eq:nonormirz}) has on integration over the moduli space of hyperbolic surfaces.   
On a hyperbolic surface, the length of an isotopy class of simple closed curves is the minimal length of a curve in that isotopy class, or equivalently the length of the unique closed geodesic in that class.  On the moduli space of hyperbolic surfaces of a given topological type the length of an isotopy class of closed curves $l([\gamma])$ is a locally defined function.  It fails to be globally defined because a curve is only well-defined up to its different images under the mapping class group, and in general $l(h\cdot[\gamma])\neq l([\gamma])$ for $h$ an element of the mapping class group.  Nevertheless, for any $f:\br^+\to\br$ one can take a sum over orbits of the mapping class group
\[F=\sum_{h\in mcg}f\left(l_{h\cdot[\gamma]}\right)\]
which is a well-defined function (when it converges) on the moduli space.  More generally, the summand can include lengths of more than one isotopy class of curves.  The McShane identities in Theorem~\ref{th:kident}, \ref{th:main} and \ref{th:kident2} are of this form, expressing a function on the moduli space as a sum of functions of lengths over orbits of the mapping class group.

In \cite{MirSim} Mirzakhani used a McShane identity on orientable hyperbolic surfaces, which is generalised by Theorem~\ref{th:main},  to calculate the volume of the moduli space of oriented hyperbolic surfaces of a given topology.  This is an application of a more general integration technique developed in \cite{MirSim}.  The moduli spaces of non-orientable surfaces have infinite volume so we integrate functions with appropriate decay.   The simplest example of such an integration over the moduli space of punctured Klein bottles is given in Section~\ref{sec:mod}.

\section{Geodesics perpendicular to the boundary}   \label{sec:proof}
In this section we prove the main theorems which are identities generalising those of  McShane and Mirzakhani.  McShane \cite{McSRem} proved his identity by considering geodesics emanating from the cusp of a hyperbolic surface.  Mirzakhani generalised this \cite{MirSim}, showing that one needs to understand geodesics perpendicular to a specified boundary component.  We will study such geodesics on any hyperbolic surface with boundary.  In particular, when the hyperbolic surface is a pair of pants or a M\"obius strip minus a disk, the functions $R$, $D$ and $E$ used in the statement of Theorem~\ref{th:main} are defined in terms of probabilities measuring the different behaviours of geodesics perpendicular to a specified boundary component.

Label the geodesic boundary components of a given hyperbolic surface $\Sigma$ by $\beta_1,...,\beta_n$, and assume $n\geq 1$.  
Take any point $x$ on the boundary component $\beta_1$ and consider the geodesic at $x$ in the direction perpendicular to $\beta_1$.  Travel along the geodesic and stop when one of the following occurs:
\begin{itemize}
\item[{\bf A}] The geodesic meets itself or $\beta_1$: stop.
\item[{\bf B}] The geodesic meets another boundary component $\beta_j$, $j>1$: stop.
\item[{\bf C}] The geodesic remains embedded for all time: don't stop.
\end{itemize}
The behaviour {\bf A}, {\bf B} and {\bf C} on any hyperbolic surface partitions $\beta_1$ into three measurable subsets.  More specifically, we will see that types {\bf A} and {\bf B} are open subsets of $\beta_1$, and quite importantly that type {\bf C} has measure zero.

\begin{figure}[ht]  
	\centerline{\includegraphics[height=2.5cm]{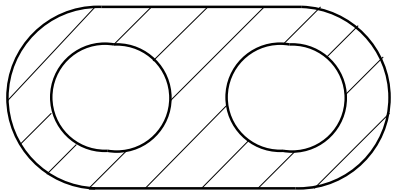}}
	\caption{Pair of pants}
	\label{fig:pants}
\end{figure}

When $\Sigma$ is a pair of pants, see Figure~\ref{fig:pants}, the boundary components are the only three embedded closed geodesics and their lengths can take on any non-negative values and uniquely determine $\Sigma$.

The probability that a geodesic perpendicular to $\beta_1$ behaves as in {\bf A} on a pair of pants
\[ \widetilde{D}(l_{\beta_1},l_{\beta_2},l_{\beta_3}):=Pr_{{\rm geodesic\ of\ type\ {\bf A}}}\Big(\pants^{\beta_2}_{\beta_3}\Big)\]
defines a function of non-negative variables since the geodesic boundary components $\beta_i$, $i=1,2,3$ can take any lengths.
Similarly define a function by the probability that a geodesic perpendicular to $\beta_1$ behaves as in {\bf B} on a pair of pants
\[ \widetilde{R}(l_{\beta_1},l_{\beta_2},l_{\beta_3}):=Pr_{{\rm geodesic\ of\ type\ {\bf B}\ meets\ }\beta_3}\Big(\pants^{\beta_2}_{\beta_3}\Big).\]

\begin{figure}[ht]  
	\centerline{\includegraphics[height=2.5cm]{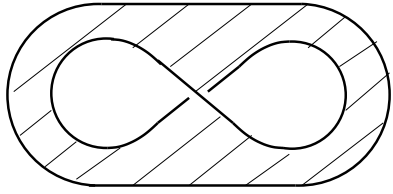}}
	\caption{M\"obius strip minus a disk}
	\label{fig:mob}
\end{figure}

When $\Sigma$ is a M\"obius strip minus a disk as shown in Figure~\ref{fig:mob}, there are exactly four embedded closed geodesics.  The lengths of any three of these geodesics determine $\Sigma$.
Denote the boundary geodesics of a M\"obius strip minus a disk by $\beta_1$ and $\beta_2$ and the interior embedded closed geodesics by $\mu$ and $\mu'$.  A trace identity, (\ref{eq:mobid}) in Section~\ref{sec:tracid} gives the relation
\begin{equation}   \label{eq:mobident} 
\cosh\left(\frac{l_{\beta_1}}{2}\right)+\cosh\left(\frac{l_{\beta_2}}{2}\right)=2\sinh\left(\frac{l_{\mu}}{2}\right)\sinh\left(\frac{l_{\mu'}}{2}\right)\end{equation}
so any three lengths determine the fourth.

A geodesic perpendicular to $\beta_1$ on a M\"obius strip minus a disk that behaves as in {\bf A} necessarily meets at least one of the one-sided embedded geodesics $\mu$ and $\mu'$.  The probability that a geodesic of type {\bf A} meets $\mu$ {\em before} it meets $\mu'$, defines the function 
\[ \widetilde{E}(l_{\beta_1},l_{\beta_2},l_{\mu}):=Pr_{{\rm geodesic\ of\ type\ {\bf A}\ meets\ }\mu{\rm\ first}}\Big(\mobius^{\beta_2}_{\mu}\Big)\]
where the diagram signifies the boundary components $\beta_1$ and $\beta_2$, and an interior embedded geodesic $\mu$.  

\subsection{Relations between probability functions}
As stated above, the behaviour {\bf A}, {\bf B} and {\bf C} of geodesics perpendicular to $\beta_1$ partitions $\beta_1$ into three subsets.  
Applying this to a pair of pants with boundary lengths $x,y,z$ with $x=l_{\beta_1}$ gives the relationship
\begin{equation} \label{eq:relDR}
\widetilde{D}(x,y,z)+\widetilde{R}(x,y,z)+\widetilde{R}(x,z,y)=1
\end{equation}
since a perpendicular geodesic is of type {\bf A} or {\bf B}, or with probability zero one of the four geodesics of type {\bf C}.

A M\"{o}bius strip minus a disk gives a relation that does not quite determine $\widetilde{E}(x,y,z)$ but gives enough information for the statement of Theorem~\ref{th:main}.
\begin{prop}
\begin{equation} \label{eq:relDER}
\widetilde{E}(x,y,z)+\widetilde{E}(x,y,z')+\widetilde{R}(x,2z,y)+\widetilde{R}(x,2z',y)=1
\end{equation}
\[{\rm for}\quad \cosh(x/2)+\cosh(y/2)=2\sinh(z/2)\sinh(z'/2).\]
\end{prop}
\begin{proof}
Denote the M\"{o}bius strip minus a disk by $\Sigma$, its boundary components by $\beta_1$, $\beta_2$ and its two interior embedded closed geodesics by $\mu$ and $\mu'$.  Put $l_{\beta_1}=x$, $l_{\beta_2}=y$, $l_{\mu}=z$ and $l_{\mu'}=z'$.  The proof uses the following two key facts.

(i) A geodesic perpendicular to $\beta_1$ of type {\bf A} meets at least one of $\mu$, $\mu'$.

(ii) A geodesic perpendicular to $\beta_1$ of type {\bf B} avoids one of $\mu$, $\mu'$.

Suppose the contrary to (i).  Then a geodesic perpendicular to $\beta_1$ meeting itself or $\beta_1$ avoids $\mu$, $\mu'$ and $\beta_2$ and hence is contained insided an annular neighbourhood of $\beta_1$.  But a simple application of the Gauss-Bonnet formula contradicts the existence of a geodesic of type {\bf A} in the annulus.  To prove (ii), cut along any embedded geodesic $\xi$ perpendicular to $\beta_1$ that meets $\beta_2$.  The remaining hyperbolic surface has the topology of a M\"{o}bius strip with convex boundary.  Since $\pi_1$ of the M\"{o}bius strip is non-trivial it contains a non-trivial embedded curve (i.e. it does not bound a disk) and since the boundary is convex there is an embedded geodesic $\nu$ in its isotopy class, so $\xi$ avoids $\nu$.  But $\Sigma$ only contains the embedded closed geodesics $\mu$ and $\mu'$ hence $\nu$ is $\mu$ or $\mu'$.

A genus zero hyperbolic surface with four boundary components is the double cover of a M\"obius strip minus a disk.  In Figure~\ref{fig:pantsmir}, $l_{\beta_1}=x$, $l_{\beta_2}=y$, and $2z$ and $2z'$ are the lengths of the double covers of the two interior geodesics which by abuse of notation we also call $\mu$ and $\mu'$.  They satisfy the hyperbolic trigonometric identity (\ref{eq:mobident}).
\begin{figure}[ht]  
	\centerline{\includegraphics[height=4cm]{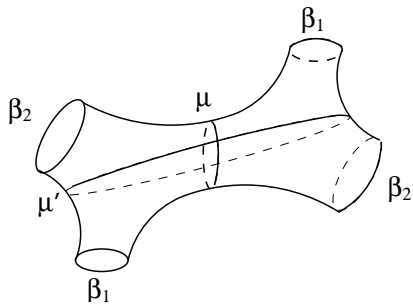}}
	\caption{Double cover of a M\"obius strip minus a disk}
	\label{fig:pantsmir}
\end{figure}

A geodesic perpendicular to $\beta_1$ is of type {\bf A} with the neighbourhood of its union with $\beta_1$ homeomorphic to a pair of pants with probability $\widetilde{E}(x,y,z)+\widetilde{E}(x,y,z')$ since it must meet $\mu$ or $\mu'$ by (i) above.  A geodesic perpendicular to $\beta_1$ is of type {\bf B}, hence stops at $\beta_2$, with probability $\widetilde{R}(x,2z,y)+\widetilde{R}(x,2z',y)$ since it must remain inside one of the two pairs of pants in the picture by (ii) above.  
With probability zero a geodesic perpendicular to $\beta_1$ is of type {\bf C} (there are eight of them.)  The proposition follows.
\end{proof}

\subsection{Proofs}  Theorems~\ref{th:main} and \ref{th:kident2} amount to understanding all types of behaviour of geodesics perpendicular to $\beta_1$ and the probability of each.  The theorems need separate treatments simply because a hyperbolic Klein bottle minus a disk is an exception to most of the topological parts of the argument for the general case.

\begin{proof}[Proof of Theorem~\ref{th:main}]

The partition of $\beta_1$ on a general hyperbolic surface into points lying on geodesics of types {\bf A}, {\bf B} and {\bf C} is better understood if we regroup them according to pairs of pants or M\"{o}bius strips minus a disk containing the geodesic.  In the special case of a Klein bottle minus a disk we will need a different grouping.

The reorganisation of probabilities goes as follows.  A pair of pants with geodesic boundary $\beta_1$, $\beta_j$ ($j>1$) and $\gamma\subset\Sigma$ of lengths $l_{\beta_1}=x$, $l_{\beta_2}=y$, $l_{\gamma}=z$, contains a geodesic perpendicular to $\beta_1$ that avoids $\gamma$ with probability
\[ \frac{R(x,y,z)}{x}:=Pr_{{\rm avoid\ }\gamma}\Big(\pants^{\beta_2}_{\gamma}\Big)=\tilde{D}(x,y,z)+\tilde{R}(x,z,y).\]
A pair of pants with geodesic boundary $\beta_1$ and $\gamma_i$, $i=1,2$, ($\gamma_i\neq\beta_j$) of lengths $l_{\beta_1}=x$, $l_{\gamma_1}=y$, $l_{\gamma_2}=z$, contains a geodesic perpendicular to $\beta_1$ that avoids $\gamma_i$, $i=1,2$, with probability
\[ \frac{D(x,y,z)}{x}:=Pr_{{\rm avoid\ }\gamma_1,\gamma_2}\Big(\pants^{\gamma_1}_{\gamma_2}\Big)=\tilde{D}(x,y,z).\]
A M\"{o}bius strip minus a disk with geodesic boundary $\beta_1$ and $\nu$, of lengths $l_{\beta_1}=x$, $l_{\nu}=y$, and containing one-sided geodesics $\mu$ and $\mu'$ of lengths $z$ and $z'$, contains a geodesic perpendicular to $\beta_1$ that avoids $\nu$, with probability
\begin{eqnarray*} 
\frac{E(x,y,z)+E(x,y,z')}{x}&:=&Pr_{{\rm avoid\ }\nu}\Big(\mobius^{\nu}_{\mu}\Big)\\
&=&\tilde{E}(x,y,z)+\tilde{E}(x,y,z').
\end{eqnarray*}
Note that this only defines the sum $E(x,y,z)+E(x,y,z')$.

On $\Sigma$, a geodesic $\xi$ perpendicular to $\beta_1$ (that stops) lies inside either a unique hyperbolic pair of pants or a unique M\"{o}bius strip minus a disk except when $\Sigma$ is a Klein bottle minus a disk.  To see this, if $\xi$ is of type {\bf A}, i.e. it stops when it meets itself or $\beta_1$, the boundary of a neighbourhood of $\xi\cup\beta_1$ consists of one or two curves.  In the former case the boundary curve is isotopic to a closed embedded geodesic in $\Sigma$ which bounds a M\"{o}bius strip minus a disk with $\beta_1$ (except when $\Sigma$ is a Klein bottle minus a disk.)  In the latter case, either the two boundary curves are isotopic to closed embedded geodesics in $\Sigma$ which bounds a pair of pants with $\beta_1$; or one of the boundary curves is isotopic to a closed embedded geodesic which bounds a M\"{o}bius strip minus a disk with $\beta_1$ (or in the special case of a Klein bottle minus a disk, neither boundary curve is isotopic to a closed embedded geodesic.)  If the geodesic $\xi$ is of type {\bf B} and meets the boundary at $\beta_j$, $j>1$, then the boundary of a neighbourhood of $\beta_1\cup\xi\cup\beta_j$ consists of a closed embedded geodesic in $\Sigma$ which bounds a pair of pants with $\beta_1$ and $\beta_j$ (except when $\Sigma$ is a M\"{o}bius strip minus a disk.)

Thus we take any pair of pants or M\"{o}bius strip minus a disk in $\Sigma$ with $\beta_1$ as a boundary component and group together all those geodesics perpendicular to $\beta_1$ of types {\bf A} and {\bf B} to get the following:
\[ 1=\sum_{\gamma_1,\gamma_2}Pr_{\rm{\bf A}}\Big(\pants^{\gamma_1}_{\gamma_2}\Big)+\sum_{j=2}^n\sum_{\gamma}Pr_{\rm{\bf B}}\Big(\pants^{\beta_j}_{\gamma}\Big)+\sum_{\nu,\mu}Pr_{\rm{\bf A}}\Big(\mobius^{\nu}_{\mu}\Big)+Pr_{\infty}(\Sigma)\]
where summands denote the probability that a geodesic perpendicular to $\beta_1$ in the pictured pair of pants or M\"{o}bius strip minus a disk is of type {\bf A} or type {\bf B}, and $Pr_{\infty}(\Sigma)$ is the probability that a geodesic is of type {\bf C}, i.e. it never stops.  For the sum over M\"{o}bius strips minus a disk, only one of the pairs $(\nu,\mu)$ and $(\nu,\mu')$ need appear (where $l_{\mu'}$ is uniquely determined by $l_{\beta_1}$, $l_{\nu}$ and $l_{\mu'}$.)  In the next expression both pairs $(\nu,\mu)$ and $(\nu,\mu')$ appear in the sum.  We have
\[  L_1=\sum_{\gamma_1,\gamma_2}D(L_1,l_{\gamma_1},l_{\gamma_2})+\sum_{j=2}^n\sum_{\gamma}R(L_1,L_j,l_{\gamma})
+\sum_{\nu,\mu}E(L_1,l_{\nu},l_{\mu})+L_1Pr_{\infty}(\Sigma)\]
where the sum is over pairs of geodesics $\gamma_1$ and $\gamma_2$ that bound a pair of pants with $\beta_1$, boundary components $\beta_j$, $j=2,..,n$ and geodesics $\gamma$ that bound a pair of pants with $\beta_1$, and one-sided geodesics $\mu$ and two-sided geodesics $\nu$ that, with $\beta_1$, bound a M\"{o}bius strip minus a disk containing $\mu$.

It remains to show that $Pr_{\infty}(\Sigma)=0$.  For a pair of pants and a M\"{o}bius strip minus a disk there are only finitely many geodesics perpendicular to $\beta_1$ that remain embedded for all time.  In general, the set of such geodesics is not even countable.  That this set has measure zero is a major part of the proof of Mirzakhani \cite{MirSim} and McShane \cite{McSRem}.  In our case, although we cannot deduce the non-orientable identity from the orientable one, the zero measure property of this set does follow from the orientable case.  Take an orientable double cover of a non-orientable hyperbolic surface so that $\beta_1$ has two preimages upstairs, one labeled $\tilde{\beta}_1$.  Geodesics perpendicular to $\beta_1$ lift to geodesics perpendicular to $\tilde{\beta}_1$.  The set of points on $\beta_1$ that lie on geodesics perpendicular to $\beta_1$ that remain embedded for all time is a subset of those points on $\tilde{\beta}_1$ that lie on geodesics perpendicular to $\tilde{\beta}_1$ that remain embedded for all time and hence the set downstairs has measure zero, i.e. $Pr_{\infty}(\Sigma)=0$.  (The sets are not necessarily equal, because a geodesic upstairs may remain embedded while its image downstairs may have self-intersection and hence should have stopped.)

Relations (\ref{eq:relDR}) and (\ref{eq:relDER}) yield the relations
$D(x,y,z)=R(x,y,z)+R(x,z,y)-x$ and $E(x,y,z)+E(x,y,z')=R(x,2z,y)+R(x,2z,y)-x$.  We arbitrarily choose $E(x,y,z)=R(x,2z,y)-x/2$ and note that this does not correspond to the probability of reaching a one-sided geodesic first.  The explicit formula 
\[ R(x,y,z)=x-\ln\frac{\cosh\frac{y}{2}+\cosh\frac{x+z}{2}}{\cosh\frac{y}{2}+\cosh\frac{x-z}{2}}\]
comes from Mirzakhani \cite{MirSim}.  It is calculated using hyperbolic trigonometry.  We will not repeat the calculation here.
\end{proof}

\begin{proof}[Proof of Theorem~\ref{th:kident2}]
For the Klein bottle minus a disk, $K$, the idea of the previous proof remains the same.  In this case all geodesics perpendicular to $\beta_1=\partial K$ are of type {\bf A}  meaning they stop because they intersect themselves or $\beta_1$.   There are three types of behaviour of these geodesics -  the geodesic intersects the unique embedded two-sided geodesic $\gamma$, zero, one or two times.

We end up with the sum of the probabilities
\[ 1=\sum_{\gamma_1,\gamma_2}Pr_{\gamma_i}\Big(\pants^{\gamma_1}_{\gamma_2}\Big)+Pr_{{\rm geodesic\ of\ type\ {\bf A}\ meets\ }\gamma}+Pr_{\infty}(\Sigma).\]
The first sum consists of a single term $\gamma_1=\gamma_2=\gamma$ and corresponds to the geodesics perpendicular to $\beta_1$ that do not meet $\gamma$.  The third term vanishes since it is the measure of a set of measure zero for the same reasons as in the previous proof.  The second sum consists of the geodesics perpendicular to $\beta_1$ that intersect $\gamma$ once or twice.  It involves a sum involving lengths of pairs of one-sided geodesics.  The summand $E(L_1,l_{\nu},l_{\mu})$ from the general case overcounts the probability that the geodesics meet $\gamma$ twice so we need to reorganise the information.  

Cut $K$ along any two disjoint one-sided geodesics $\mu$ and $\mu'$ to get a pair of pants with boundary lengths $x=l_{\beta_1}$, $2y=2l_{\mu}$ and $2z=2l_{\mu'}$.  The two closest points on the boundary component of length $2y$ divide it into two equal parts, and in Figure~\ref{fig:pantsf} we divide these again to get four parts of equal length $y/2$.  Similarly do this on the other boundary component.
\begin{figure}[ht]  
	\centerline{\includegraphics[height=4cm]{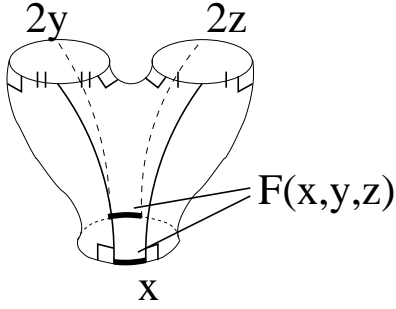}}
	\caption{Probability on the Klein bottle minus a disk}
	\label{fig:pantsf}
\end{figure}

The shaded region on $\beta_1$ of total length $F(x,y,z)$ denotes the probability that a geodesic perpendicular to $\beta_1$ on $K$ will avoid at least one of $\mu$ and $\mu'$, and if it meets $\mu'$, say, then it has already crossed $\mu''$, the other one-sided geodesic disjoint from $\mu$.  This last requirement simply means that both $F(L_1,\mu,\mu')$ and $F(L_1,\mu,\mu'')$ appear in the sum of probabilities since a perpendicular geodesic meets one of $\mu$, $\mu'$ first.

Thus we have
\[  L_1=D(L_1,l_{\gamma},l_{\gamma})+\sum_{\mu,\mu'}F(L_1,l_{\mu},l_{\mu'})\]
where the sum is over all pairs of disjoint one-sided geodesics $\mu$ and $\mu'$.  To calculate $F(x,y,z)$ we use hyperbolic trigonometry.  
\begin{figure}[ht]  
	\centerline{\includegraphics[height=4cm]{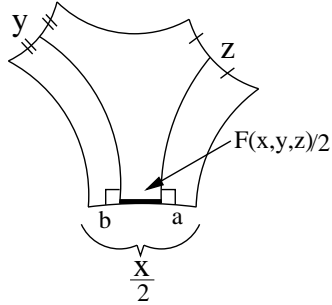}}
	\caption{Right-angle hexagon}
	\label{fig:hypf}
\end{figure}

Figure~\ref{fig:hypf} shows a hyperbolic right-angled hexagon with perpendiculars dropped from the midpoints of two sides to the third side.  The perpendiculars remain inside the hexagon and do not meet each other by a simple application of the Gauss-Bonnet formula.  To calculate $F(x,y,z)$ we calculate $a$ and $b$ and use $a+F(x,y,z)/2+b=x/2$.  The two quadrilaterals containing the lengths $a$ and $b$ have 3 right-angles and hence satisfy a simple hyperbolic trigonometric identity:
\[ \tanh{a}=\tanh(z/2)/\cosh{y'},\quad \tanh{b}=\tanh(y/2)/\cosh{z'}\] 
for
\[ \cosh{y'}=\frac{\cosh{y}+\cosh(x/2)\cosh{z}}{\sinh(x/2)\sinh{z}},\quad
\cosh{z'}=\frac{\cosh{z}+\cosh(x/2)\cosh{y}}{\sinh(x/2)\sinh{y}}\]
where $y'$ and $z'$ are the lengths of sides of the hexagon opposite $y$ and $z$, respectively, and satisfy a hyperbolic trigonometric identity on the hexagon.  So
\[ \tanh{a}=\frac{\tanh(z/2)\sinh(x/2)\sinh{z}}{\cosh{y}+\cosh(x/2)\cosh{z}},\quad \tanh{b}=\frac{\tanh(y/2)\sinh(x/2)\sinh{y}}{\cosh{z}+\cosh(x/2)\cosh{y}}\] 
and
\begin{eqnarray*} 
F(x,y,z)&=&x-2a-2b\\
&=&\frac{x}{2}-\ln\frac{\cosh{y}+\exp\frac{x}{2}\cosh{z}-\sinh\frac{x}{2}}{\cosh{y}+\exp\frac{-x}{2}\cosh{z}+\sinh\frac{x}{2}}.
\end{eqnarray*}
Since 
\[ L-D(L,l_{\gamma},l_{\gamma})=2\ln\frac{1+e^{L/2}e^{l_{\gamma}}}{e^{L/2}+e^{l_{\gamma}}}\] 
the sum of probabilities becomes
\[\sum_{\mu,\nu}F(L,l_{\mu},l_{\nu})=2\ln\frac{1+e^{L/2}e^{l_{\gamma}}}{e^{L/2}+e^{l_{\gamma}}}\]
and the theorem is proven.
\end{proof}

\section{Trace identities}  \label{sec:tracid}

Lengths of closed geodesics on a hyperbolic surface can be expressed in terms of traces of elements in $PGL(2,\br)$.  Trace identities between elements in $PGL(2,\br)$ then give identities between lengths of closed geodesics.  In this section we prove trace identities required throughout the paper, and mainly focus on the less familiar non-orientable cases. 

A closed geodesic $\gamma$ on a hyperbolic surface $\Sigma$ corresponds to a conjugacy class in $\pi_1\Sigma$ and hence it is sent to a conjugacy class represented by $A_{\gamma}\in PGL(2,\br)$, under the representation $\pi_1\Sigma\to PGL(2,\br)$ that defines the hyperbolic structure.  On an orientable hyperbolic surface the representation takes its values in $PGL^+(2,\br)$.  If we take a matrix representative $A_{\gamma}$ satisfying $\det A_{\gamma}=1$ and $\tr A_{\gamma}>0$, there is a relationship
\begin{equation}  \label{eq:trorient}
 \tr A_{\gamma}=2\cosh\frac{1}{2}l_{\gamma}.
\end{equation}
To see this, conjugate $A_{\gamma}$ so that it preserves the imaginary line in upper half space $\bh^+$.  Then
\[ A_{\gamma}=\left(\begin{array}{cc}a&0\\0&1/a\end{array}\right)\]
so $A_{\gamma}\cdot i=a^2i$ which is translation by the hyperbolic length $\ln a^2$.  In other words $l_{\gamma}=\ln a^2$ and (\ref{eq:trorient}) follows.

On a non-orientable surface (\ref{eq:trorient}) is still true for two-sided closed geodesics.   For one-sided closed geodesics we have
\begin{equation}  \label{eq:trnonorient}
 \tr A_{\gamma}=2\sinh\frac{1}{2}l_{\gamma}.
\end{equation}
where we have scaled $A_{\gamma}$ to have determinant $-1$ and positive trace.  The change from $\cosh$ to $\sinh$ is due to the fact that a one-sided homotopy class maps to a negative determinant matrix $A$, or equivalently it has non-trivial image under the homomorphism 
\[ w_1:\pi_1\Sigma\to PGL(2,\br)\stackrel{sgn\det}{\longrightarrow}\bz_2.\]
A negative determinant matrix $A$ acts on $\bh^+$ by 
\[ z\mapsto A\bar{z}\]
and as above, we arrange $A_{\gamma}$ so that it preserves the imaginary line in upper half space $\bh^+$ to get
\[ A_{\gamma}=\left(\begin{array}{cc}a&0\\0&-1/a\end{array}\right)\]
so $A_{\gamma}\cdot i=a^2i$ and as before $l_{\gamma}=\ln a^2$ yielding (\ref{eq:trnonorient}).

Relations (\ref{eq:trorient}) and (\ref{eq:trnonorient}) are most naturally expressed as
\[ \frac{(\tr A_{\gamma})^2}{\det A_{\gamma}}=2\pm2\cosh{l_{\gamma}}\] 
since this is well-defined on $PGL(2,\br)$ (where $\pm=sgn\det(A_{\gamma})$.)  We usually rescale so that $\det A_{\gamma}=\pm1$ for ease of use.  

The following trace identity is better known on $PSL(2,\br)$ but we write it on $PGL(2,\br)$ for naturality and generality, before specialising by rescaling.
For $A,B\in PGL(2,\br)$, and $B^{\dagger}$, the adjugate of $B$, defined by $BB^{\dagger}=\det{B}$ 
\begin{equation}  \label{eq:tracident}
\det B\cdot(\tr A)^2+\det A\cdot(\tr B)^2+(\tr AB)^2-\tr A\cdot \tr B\cdot \tr AB=\tr ABA^{\dagger}B^{\dagger}+2\det A\det B.
\end{equation}
(One can prove (\ref{eq:tracident}) by extending $\dagger$ to a linear map on $M(2,\br)$ and converting the quadratic identity to a bilinear identity, or by directly verifying it with the computer.)  

The trace identity (\ref{eq:tracident}) applied to $A=A_{\gamma}$ and $B=A_{\mu}$ has different geometric interpretations according to the behaviour of the closed geodesics $\gamma$ and $\mu$.  It neatly encodes some difficult hyperbolic trigonometric identities.  We will apply it to simple closed curves intersecting exactly once.

Let $\gamma$ and $\mu$ be two-sided simple closed geodesics that intersect exactly once (and hence both are non-separating.)  Then the commutator $\big[[\gamma],[\mu]\big]\in\pi_1\Sigma$ is represented by a two-sided simple closed geodesic that bounds a neighbourhood of $\gamma\cup\mu$ which is homeomorphic to $T^2-D^2$.   The hyperbolic structure gives a local representation $\pi_1\Sigma\to PSL(2,\br)$ (since two-sided curves map to $PGL^+(2,\br)\subset PGL(2,\br)$.)  In this case (\ref{eq:tracident}) becomes the better known
\[ (\tr A_{\gamma})^2+(\tr A_{\mu})^2+(\tr A_{\gamma}A_{\mu})^2-\tr A_{\gamma}\cdot \tr A_{\mu}\cdot \tr A_{\gamma}A_{\mu}=\tr [A_{\gamma},A_{\mu}]+2.\]
The sign of $\tr[A_{\gamma},A_{\mu}]$ is well-defined on 
$PSL(2,\br)$ and in fact $\tr[A_{\gamma},A_{\mu}]<0$ which uses the following continuity argument.  The space of geometric pairs $(A_{\gamma},A_{\mu})$ is connected so since $\tr[A_{\gamma},A_{\mu}]$ remains away from zero (the representation is geometric (Fuchsian) meaning $|\tr[A_{\gamma},A_{\mu}]|\geq 2$) it is sufficient to calculate the sign of $\tr[A_{\gamma},A_{\mu}]$ in one example (which we do not do here, although we will construct an example for the non-orientable case.)  Using the quadratic identity, $\tr[A_{\gamma},A_{\mu}]\leq-2$ implies that $\tr A_{\gamma}\cdot \tr A_{\mu}\cdot \tr A_{\gamma}A_{\mu}>0$.  

Lift the local representation $\pi_1\Sigma\to PSL(2,\br)$ to $SL(2,\br)$ so that $\tr A_{\gamma}>0$ and $\tr A_{\mu}>0$.  From $\tr A_{\gamma}\cdot \tr A_{\mu}\cdot \tr A_{\gamma}A_{\mu}>0$ this forces $\tr A_{\gamma}A_{\mu}>0$ and we have
\begin{equation}  \label{eq:markoff} 
x^2+y^2+z^2-xyz=-\partial+2,\quad\left\{\begin{array}[m]{ll}x=2\cosh\frac{1}{2}l_{\gamma},&  z=2\cosh(l_{[\gamma][\mu]}/2),\\y=2\cosh(l_{\mu}/2),& \partial=2\cosh(l_{[[\gamma],[\mu]]}/2).
\end{array}\right.
\end{equation}
When $\delta=2$, so the right hand side of (\ref{eq:markoff}) is zero, such a triple $(x,y,z)$ is known as a {\em Markoff triple.}  

Now let $\gamma$ be a one-sided simple closed geodesic and $\mu$ a two-sided simple closed geodesic that intersect exactly once.  Then the commutator $\big[[\gamma],[\mu]\big]\in\pi_1\Sigma$ no longer has good geometric meaning, instead the product $[\gamma][\mu][\gamma]^{-1}[\mu]$
is represented by a two-sided simple closed geodesic that bounds a neighbourhood of $\gamma\cup\mu$ which is homeomorphic to the Klein bottle minus a disk $K-D^2$.  Apply the identity 
\begin{equation}  \label{eq:mobid}
\tr\alpha\beta+\tr\alpha\beta^{\dagger}=\tr\alpha\cdot\tr\beta
\end{equation} 
to $\alpha=ABA^{\dagger}$ and $\beta=B$ to get
\[\tr ABA^{\dagger}B+\tr ABA^{\dagger}B^{\dagger}=
\tr ABA^{\dagger}\cdot\tr B=\det A\cdot(\tr B)^2\]
so (\ref{eq:tracident}) becomes
\begin{equation}  \label{eq:tracident2}
\det B\cdot(\tr A)^2+(\tr AB)^2-\tr A\cdot \tr B\cdot \tr AB=-\tr ABA^{\dagger}B+2\det A\det B
\end{equation}
which has replaced $\tr ABA^{\dagger}B^{\dagger}$ by $\tr ABA^{\dagger}B$, the former geometrically meaningful in the orientable case, and the latter geometrically meaningful in the nonorientable case.

As in the orientable case, the sign of $\tr ABA^{\dagger}B$ is well-defined on $PGL(2,\br)$ and by a continuity argument on the connected space of geometric pairs $(A,B)$ (i.e. pairs corresponding to a one-sided geodesic intersecting a two-sided geodesic once) it is determined by calculating a single example.  Such an example is given by the subgroup of $PGL(2,\bz)$ shown in (\ref{eq:intrep}).  In this example, $\tr ABA^{\dagger}B>0$ thus this is true for all geometric pairs $(A,B)$.  The quadratic identity then implies $\tr A\cdot\tr B\cdot\tr AB>0$.

We lift the representation $\pi_1\Sigma\to PGL(2,\br)$ to $GL(2,\br)$ and rescale so that $\det A_{\gamma}=-1$, $\det A_{\mu}=1$, $\tr A_{\gamma}>0$ and $\tr A_{\mu}>0$.  This forces $\tr A_{\gamma}A_{\mu}>0$ so (\ref{eq:tracident2}) becomes
\[
Y_1^2+Y_2^2-Y_1Y_2Z=-\partial-2,\quad
\left\{\begin{array}[m]{ll}Y_1=2\sinh\frac{1}{2}l_{\gamma},&Y_2=2\sinh(l_{[\gamma][\mu]}/2),\\ Z=2\cosh(l_{\mu}/2),& \partial=2\cosh(l_{[\gamma][\mu][\gamma]^{-1}[\mu]}/2).\end{array}\right.
\]

\begin{figure}[ht]  
	\centerline{\includegraphics[height=4cm]{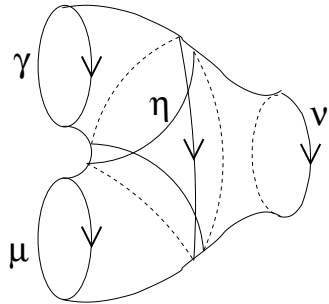}}
	\caption{$l_{\eta}=S(l_{\gamma},l_{\mu},l_{\nu})$}
	\label{fig:3fold}
\end{figure}

A third geometric case when $\gamma$ and $\mu$ are disjoint two-sided simple closed geod\-esics is not important for the paper but we include it out of interest.  In this case $A_{\gamma}A_{\mu}=A_{\nu}$ for $\nu$ a simple closed geodesic that bounds a pair of pants with $\gamma$ and $\mu$, and the trace identity (\ref{eq:tracident}) expresses the length of the closed geodesic with (topological) three-fold symmetry, drawn as $\eta$ in Figure~\ref{fig:3fold}, as a symmetric function in the boundary lengths.

\section{Punctured Klein bottle}  \label{sec:klein}
Theorem~\ref{th:kident} follows from Theorem~\ref{th:kident2} if we divide both sides by $L$ and take the limit $L\to 0$.  
In this section we give an alternative, elementary proof of Theorem~\ref{th:kident} similar to Bowditch's proof \cite{BowPro} of McShane's original identity on the punctured torus.  As with a punctured torus, the lengths of all geodesics on a hyperbolic punctured Klein bottle can be calculated recursively (and hence easily listed on the computer.)

Keeping the notation of the theorem, $\gamma$ denotes the unique two-sided embedded closed geodesic, and $\gamma_i$ denotes the sequence of one-sided embedded closed geodesics in the hyperbolic punctured Klein bottle.  Put 
\[ Z=2\cosh\frac{1}{2}l_{\gamma},\quad  y_i=\sinh\frac{1}{2}l_{\gamma_i}.\]

The trace identities give the following relations.  The puncture has zero length so the quadratic identity becomes
\begin{equation}  \label{eq:quad}
y_i^2+y_{i+1}^2-y_iy_{i+1}Z=-1.
\end{equation}
We can replace $y_{i+1}$ in this quadratic identity by $y_{i-1}$ and hence they are the two roots of a common quadratic equation.  Thus
\begin{equation}  \label{eq:quad1}
y_{i-1}+y_{i+1}=y_iZ\quad{\rm and}\quad y_{i-1}y_{i+1}-y_i^2=1.
\end{equation}
Use these to evaluate the series
\begin{eqnarray*}
\sum_{i=-\infty}^{\infty}\frac{1}{1+y_i^2+y_{i+1}^2}&=&
\sum_{i=-\infty}^{\infty}\frac{1}{y_iy_{i+1}Z}\\
&=&\frac{1}{Z}\sum_{i=-\infty}^{\infty}\left(\frac{y_{i-1}}{y_i}-\frac{y_i}{y_{i+1}}\right)\\
&=&\frac{1}{Z}\left(\lim_{i\to-\infty}\frac{y_{i-1}}{y_i}-\lim_{i\to\infty}\frac{y_i}{y_{i+1}}\right).
\end{eqnarray*}
From (\ref{eq:quad}) we have 
\[ \left(\frac{y_i}{y_{i+1}}\right)^2+1-Z\frac{y_i}{y_{i+1}}=-\frac{1}{y_{i+1}^2}\]
and $1/y_i\to 0$ as $i\to\pm\infty$ (the lengths of the geodesics go to infinity) so
\[ \lim_{i\to\pm\infty}y_i/y_{i+1}=\lambda_{\pm}\] satisfies
$\lambda_{\pm}^2-Z\lambda_{\pm}+1=0\Rightarrow \lambda=\exp(\pm l_{\gamma}/2)$.

From $y_{i-1}y_{i+1}-y_i^2=1$ we calculate recursively $\displaystyle{y_{i+1}=\frac{1+y_i^2}{y_{i-1}}}$
so if $y_i>y_{i-1}$ then $y_{i+1}>y_i$.  Similarly, if $y_i>y_{i+1}$ then $y_{i-1}>y_i$.  Eventually as $i$ gets large the sequence is strictly increasing and as $-i$ becomes large the sequence is strictly decreasing (since a hyperbolic surface has a minimal length closed geodesic), so  $\lim_{i\to\infty}y_i/y_{i+1}=\exp-\frac{1}{2}l_{\gamma}$ and 
$\lim_{i\to-\infty}y_i/y_{i+1}=\exp\frac{1}{2}l_{\gamma}$.  Thus
\begin{eqnarray*}
 \frac{1}{Z}\left(\lim_{i\to-\infty}\frac{y_{i-1}}{y_i}-\lim_{i\to\infty}\frac{y_i}{y_{i+1}}\right)&=&\frac{1}{Z}\left(\exp\frac{1}{2}l_{\gamma}-\exp-\frac{1}{2}l_{\gamma}\right)\\
 &=&
\frac{2\sinh\frac{1}{2}l_{\gamma}}{2\cosh\frac{1}{2}l_{\gamma}}=\tanh\frac{1}{2}l_{\gamma}.
\end{eqnarray*}
Hence
\[ \sum_{i=-\infty}^{\infty}\frac{1}{1+y_i^2+y_{i+1}^2}=\tanh\frac{1}{2}l_{\gamma}\]
which is the identity (\ref{eq:kident}).

\subsection{Rational solutions}
Surfaces with rational traces correspond to rational solutions of (\ref{eq:quad}).  If $y_n$ and $y_{n+1}$ are rational, then all $y_i$ and $Z$ are rational.  There is exactly one hyperbolic surface with all traces integral.  This has
\[ y_i=F_{2i},\quad i\geq 0,\quad y_{-i}=F_{2i-2},\quad i>0,\]
where $F_i$ is the $i$th Fibonacci number (for $F_0=1=F_1$.)  The triples $(F_{2i},F_{2i+2},3)$ give all integral solutions to (\ref{eq:quad}).  The hyperbolic surface with integral traces is a quotient of $\bh^2$ by the index 12 subgroup of $PGL(2,\bz)$:
\begin{equation}  \label{eq:intrep}
\Gamma=\langle A,B\rangle\subset PGL(2,\bz),\quad A=\left(\begin{array}{cc}0&-1\\-1&2\end{array}\right),\quad B=\left(\begin{array}{cc}1&1\\1&2\end{array}\right).
\end{equation}  
It is analogous to the modular curve given by the quotient of $\bh^2$ by the index 6 commutator subgroup of $PSL(2,\bz)$ which is a hyperbolic punctured torus with all traces integral.

\subsection{Extended Quasi-Fuchsian 3-manifolds}

A complete hyperbolic 3-manifold $M=\bh^3/\Gamma$ is the quotient of hyperbolic space by a discrete subgroup $\Gamma\subset PSL(2,\bc)$ which is the image of a representation of $\pi_1M$. It has conformal boundary $\Omega(\Gamma)/\Gamma$ where $\Omega(\Gamma)$ is the largest open subset of $S^2_{\infty}$ (the conformal boundary of $\bh^3$) on which $\Gamma$ acts properly discontinuously.  The limit set $\Lambda(\Gamma)=S^2_{\infty}-\Omega(\Gamma)$ is the set of accumulation points of an arbitrary point in $\bh^3$.  Hyperbolic surfaces correspond to Fuchsian groups which have $\Lambda(\Gamma)$ equal to a round circle, or equivalently take values in (a conjugate of) $PSL(2,\br)$.   More generally, if $\Lambda(\Gamma)$ is a Jordan curve homeomorphic to a circle and $\Omega(\Gamma)$ consists of two disjoint sets, then either $\Gamma$ preserves the two components of $\Omega(\Gamma)$, in which case it is quasi-Fuchsian and $M$ is homeomorphic to a surface times $\br$, or $\Gamma$ swaps the two components of $\Omega(\Gamma)$, in which case it is {\em extended quasi-Fuchsian} and $M$ is homeomorphic to a non-trivial $\br$ bundle over a non-orientable surface.  Its conformal boundary is the orientable double cover of the non-orientable surface.

\begin{proof}[Proof of Theorem~\ref{th:3man}]

The identity (\ref{eq:kident}) extends to allow complex values of $y_i$.  To prove the more general version of the identity, it is sufficient to show the series on the left hand side of (\ref{eq:kident}) converges, since the proof of the identity goes through without change.  Recall that given initial values $y_0,y_1\in\bc$, they determine $l_{\gamma}\in\bc$ via (\ref{eq:quad}) and $\{y_i\}$ via either equation in (\ref{eq:quad1}).  To analyse the convergence of the series (\ref{eq:kident}) we need to understand the growth of the $|y_i|$ as $i\to\pm\infty$, and we do this by giving a general solution of $\{y_i\}$ from $y_0,y_1\in\bc\backslash\{0\}$.
\begin{equation}   \label{eq:gensol}
y_n=A\exp\frac{n}{2}l_{\gamma}+B\exp-\frac{n}{2}l_{\gamma},\quad A+B=y_0,\quad AB=4/\sinh^2\frac{1}{2}l_{\gamma}.
\end{equation}
Note that $y_1$ enters (\ref{eq:gensol}) by together with $y_0$ determining $l_{\gamma}$.  The proof of this general solution is by direct substitution into (\ref{eq:quad}) or (\ref{eq:quad1}), and the fact that there exists a unique solution given $y_0$ and $y_1$.  Now we put the following assumption on $y_0$ and $y_1$ which we interpret geometrically below.
\begin{condition}  \label{th:cond}
Restrict to pairs $y_0,y_1\in\bc\backslash\{0\}$ such that $Re(l_{\gamma})\neq 0$ for $l_{\gamma}$ defined by $\cosh\frac{1}{2}l_{\gamma}=(1+y_0^2+y_1^2)/(2y_0y_1)$.
\end{condition}
Now, put $Re(l_{\gamma})=\lambda\neq0$.
From (\ref{eq:gensol}) 
\[ |y_n|>\big||A|\exp(n\lambda/2)-|B|\exp(-n\lambda/2)\big|\sim\exp(n|\lambda|/2)\]
since both $|A|$ and $|B|$ are non-zero and one of the terms grows exponentially while the other decays exponentially.  Thus $|y_n|$ grows exponentially.  Now replace the summands in the series (\ref{eq:kident}) by the more convenient $1/y_iy_{i+1}Z$ so we see that the terms decay exponentially, the series converges absolutely, and the proof of the identity goes through under the assumption of Condition~\ref{th:cond}.

In the hyperbolic 3-manifold homeomorphic to the non-trivial $\br$ bundle over the punctured Klein bottle there is a unique closed geodesic $\gamma$ that projects to an embedded two-sided curve on the surface.  It corresponds to a loxodromic element in $PSL(2,\bc)$.  A {\em loxodromic} element $A\in PSL(2,\bc)$ satisfies $tr(A)\notin[-2,2]$.  Thus $2\cosh(l_{\gamma})\notin[-2,2]$.  But $\cosh(a+ib)=\cosh{a}\cos{b}+i\sinh{a}\sin{b}$ so $2\cosh(a+ib)\in[-2,2]$ is equivalent to $a=0$ and thus $l_{\gamma}$ corresponding to a loxodromic element is equivalent to Condition~\ref{th:cond}.  The theorem follows.
\end{proof}

Theorem~\ref{th:3man} is an analogue of Bowditch's study of two-generator subgroups of $PSL(2,\bc)$ using Markoff triples \cite{BowMar}.  In that case, it is difficult to tell when a Markoff triple---$(x,y,z)$ the norm of each greater than 2 and satisfying (\ref{eq:markoff}) for $\delta=2$---corresponds to a Kleinian group.  The situation for the punctured Klein bottle is much easier, since Condition~\ref{th:cond} gives the complete answer.

\section{Moduli space}  \label{sec:mod}

For $n+k-2>0$, denote by $\modm_{k,n,\pm}({\bf L})$ the moduli space of hyperbolic surfaces with $n$ geodesic boundary components of lengths ${\bf L}=(L_1,...,L_n)$ and fixed topology
\begin{equation}  \label{eq:top}
\Sigma\cong\Sigma'_m\#(k-2m)\br\bp^2,\quad \Sigma' 
{\rm\ orientable\ of\ genus\ }m
\end{equation} 
where $\pm$ indicates when the surface is orientable or non-orientable.  Note that $k=2m$ in the orientable case and in the non-orientable case the topology is independent of $m$.  The moduli space is an orbifold of dimension 
\[\dim\modm_{k,n,\pm}({\bf L})=3k-6+2n=-3\chi-n\]
where $\chi$ is the Euler characteristic.
The map
\begin{equation}  \label{eq:2to1}
\modm_{g,n}({\bf L})\stackrel{2:1}{\longrightarrow}\modm_{2g,n,+}({\bf L})
\end{equation}
from the usual moduli space of {\em oriented} hyperbolic surfaces is generically 2:1 since there are two orientations on each orientable surface.  (A half dimensional space of surfaces possesses an orientation preserving isometry that switches the orientations.  When ${\bf L}=0$, these are the real algebraic curves, or Klein surfaces.)  Below we give an analogue of (\ref{eq:2to1}) for non-orientable surface by equipping a surface with an orientation of one of its boundary components.

To see the manifold, or orbifold, structure of the moduli space locally, and in particular its dimension, as usual we take a maximal set of disjoint embedded closed geodesics on a given surface.  Deformations of the hyperbolic structure are obtained by varying the lengths of these geodesics.  Further deformations are obtained by cutting along two-sided geodesics, rotating and gluing back.  The gluing along a one-sided geodesic is rigid and hence does not allow rotations, i.e. a one-sided geodesic must glue to itself via the antipodal map.  A surface of type $(k,n,-)$ has the topology of (\ref{eq:top}), so it possesses a maximal set of disjoint closed embedded geodesics consisting of $k-2m$ one-sided geodesics and $m-3+n+k$  two-sided geodesics.  Thus the dimension of the deformation space is $k-2m+2(m-3+n+k)=3k-6+2n$ as given above.

If we cut along the maximal set of disjoint embedded closed geodesics we get a pair of pants decomposition of the surface.  The rotation angle along a two-sided geodesic can be given a value that is almost well-defined, by measuring the angle difference between distinguished points on each side of the geodesic, given by the two closest points to the two other boundary components of the pair of pants.  The angle $\theta({\rm mod\ }l)$ along a two-sided geodesic of length $l$ is well-defined up to
\begin{equation}  \label{eq:theta}
\theta\mapsto l-\theta,\quad\theta\mapsto \frac{l}{2}+\theta,\quad\theta\mapsto \frac{l}{2}-\theta.
\end{equation}

\subsection{Volume form}
The moduli space of hyperbolic surfaces with a given topology comes equipped with a local volume form which is globally well-defined up to sign.
\begin{theorem}  \label{th:vol}
Up to sign, the volume form
\begin{equation}  \label{eq:vol}
d{\rm vol}=\bigwedge_{\gamma\ 2-{\rm sided}}(dl_{\gamma}\wedge d\theta_{\gamma})\ \wedge\bigwedge_{\mu\ 1-{\rm sided}}\coth(l_{\mu}/2)dl_{\mu}
\end{equation}
is independent of the choice of pair of pants decomposition.
\end{theorem}
\begin{proof}
The volume form is only well-defined up to sign for two reasons.   The rotation angle $\theta_{\gamma}$ along the geodesic $\gamma$ is well-defined up to the transformations (\ref{eq:theta}) which can transform $d\theta_{\gamma}\mapsto-d\theta_{\gamma}$.  Also, the order of the $\coth(l_{\mu})dl_{\mu}$ terms is not well-defined so switching two can alter the sign.

The main point of the theorem is that up to sign the volume form is invariant under the mapping class group, and it is independent of the choice of the topological pair of pants decomposition.  (For example, the closed orientable genus 2 surface has two topologically distinct pair of pants decompositions.)  We prove both of these at the same time by showing that the volume form is invariant, up to sign, under local changes in the pair of pants decomposition.  A sequence of such changes can give any element of the mapping class group.

There are four local changes of a pair of pants decomposition.  Two cases involve only two-sided geodesics - replace a two-sided geodesic with another intersecting two-sided geodesic and consider the two cases when the geodesic is non-separating or separating.  We will not consider these two cases here since they follow from \cite{WolWei}.

The other two cases consist of 

(i) replace a one-sided geodesic $\mu$ by another one-sided geodesic $\mu'$ that intersects $\mu$ exactly once;

(ii) replace two non-intersecting one-sided geodesics by a two-sided geodesic that intersects both of the one-sided geodesics once.

For case (i), thicken $\mu\cup\mu'$ to get a M\"{o}bius strip minus a disk $\Sigma$.  The two boundary components of $\Sigma$ are isotopic to geodesics so we can choose $\Sigma$ to have geodesic boundary $\beta_1$ and $\beta_2$.  The trace identity (\ref{eq:mobident}) $\cosh\left(\frac{l_{\beta_1}}{2}\right)+\cosh\left(\frac{l_{\beta_2}}{2}\right)=2\sinh\left(\frac{l_{\mu}}{2}\right)\sinh\left(\frac{l_{\mu'}}{2}\right)$ and the fact that $\beta_1$ and $\beta_2$ are unchanged by the swap of $\mu$ with $\mu'$ we have
$d\log\sinh\left(\frac{l_{\mu}}{2}\right)+d\log\sinh\left(\frac{l_{\mu'}}{2}\right)=0$
and hence 
\begin{equation}  \label{eq:symp1}
\coth(l_{\mu})dl_{\mu}=-\coth(l_{\mu'})dl_{\mu'}.
\end{equation}  
Furthermore, not only do the lengths of $\beta_1$ and $\beta_2$ remain constant, but also the rotation angles at each of these geodesics remain constant.  This is because the distinguished points on $\beta_1$ and $\beta_2$ marked by shortest geodesics between boundary components of a pair of pants decomposition, are the same for the two pair of pants decomposition obtained by cutting along 
$\mu$ or $\mu'$.  This fact can be seen in Figure~\ref{fig:pantsmir} which shows the double cover of a M\"{o}bius strip minus a disk, and hence both pair of pants decompositions on the same diagram.  Thus all other coordinates are unchanged and the volume element changes by $\coth(l_{\mu}/2)dl_{\mu}\mapsto\coth(l_{\mu'})dl_{\mu'}=-\coth(l_{\mu})dl_{\mu}$ which is multiplication by -1.

For case (ii), take any path $p$ running between the non-intersecting one-sided geodesics $\mu_1$ and $\mu_2$, so that $p$ intersects each $\mu_i$ only at the endpoints of $p$.  Thicken $\mu_1\cup\mu_1\cup p$ to get a Klein bottle minus a disk $K$ and as before we can choose $K$ so that it has geodesic boundary.  Inside $K$ there is a unique two-sided geodesic $\gamma$.  We wish to replace $\mu_1$ and $\mu_2$ by $\gamma$.  This gives a new pair of pants decomposition with coordinates the length $l_{\gamma}$ and rotation angle $\theta_{\gamma}$ replacing the lengths $l_{\mu_1}$ and $l_{\mu_2}$.  The coordinates are shown together in the following diagram. 
\begin{figure}[ht]  
	\centerline{\includegraphics[height=4cm]{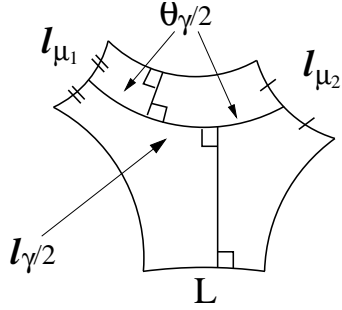}}
	\caption{Coordinates $l_{\gamma}$, $\theta_{\gamma}$, $l_{\mu_1}$ and $l_{\mu_2}$}
	\label{fig:dvol}
\end{figure}

Figure~\ref{fig:dvol} shows a right-angle hexagon, two copies giving the Klein bottle cut along the one-sided geodesics $\mu_1$ and $\mu_2$.  The geodesic joining the midpoints of the sides $\mu_1$ and $\mu_2$ gives half of the unique two-sided geodesic $\gamma$.  The length $\theta_{\gamma}/2$ appears twice in the diagram as either of the intervals of equal length.  
\begin{figure}[ht]  
	\centerline{\includegraphics[height=2cm]{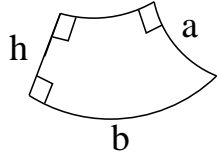}}
	\caption{$\sinh{h}\cdot\cosh{b}=\sinh{a}$}
	\label{fig:quad}
\end{figure}

Put $Y_i=\sinh(l_{\mu_i}/2)$ and apply the hyperbolic trigonometric identity in Figure~\ref{fig:quad} to the two quadrilaterals in Figure~\ref{fig:dvol} to get
\begin{equation}  \label{eq:quadrel} 
\frac{Y_1}{Y_2}=\frac{\cosh(\theta_{\gamma}/2)}{\cosh(l_{\gamma}/2-\theta_{\gamma}/2)}.
\end{equation}
Also the trace identity (\ref{eq:tracident2}) gives
\[ Y_1^2+Y_2^2-2Y_1Y_2\cosh\frac{1}{2}l_{\gamma}={\rm constant}\]
where the right hand side depends on the length of the boundary of $K$ which is constant.  So we substitute (\ref{eq:quadrel}) to get
\[\frac{\rm constant}{Y_1^2}=1+\frac{\cosh^2(l_{\gamma}/2-\theta_{\gamma}/2)}{\cosh^2(\theta_{\gamma}/2)}-\frac{2\cosh(l_{\gamma}/2-\theta_{\gamma}/2)\cosh\frac{1}{2}l_{\gamma}}{\cosh(\theta_{\gamma}/2)}\] and get rid of the constant by taking the logarithmic derivative.  Together with the logarithmic derivative of (\ref{eq:quadrel}) we get
\begin{equation} \label{eq:symp2}
dl_{\gamma}\wedge d\theta_{\gamma}=dY_1\wedge dY_2=
\coth(l_{\mu_1}/2)\coth(l_{\mu_2}/2)dl_{\mu_1}\wedge dl_{\mu_2}.
\end{equation}
Under the change from $\mu_1,\mu_2$ to $\gamma$ the length of the boundary of $K$ is unchanged but unlike the previous case, the angle associated to the boundary does change.  The change is given by one of the two subintervals of the length $L$ in Figure~\ref{fig:dvol}.  Label this change by $\phi$, then $\phi=\phi(L,l_{\mu_1},l_{\mu_2})$ so in particular $d\phi$ is a linear combination of $dl_{\mu_1}$ and $dl_{\mu_2}$, or equivalently a linear combination of $dl_{\gamma}$ and $d\theta_{\gamma}$, so the volume element is unchanged by the change $d\phi$ and hence $d{\rm vol}$ changes by (\ref{eq:symp2}) and the theorem is proven.
\end{proof}

\subsection{The oriented moduli space}

The moduli space of orientable hyperbolic surfaces $\modm_{k,n,+}({\bf L})$ possesses a double cover $\modm_{g,n}({\bf L})$ of {\em oriented} hyperbolic surfaces of genus $g=k/2$.   On $\modm_{g,n}({\bf L})$ the volume form (\ref{eq:vol}) is well-defined without the $\pm 1$ ambiguity.  Teichm\"uller space is a cover of the oriented moduli space consisting of oriented hyperbolic surfaces with a marking, i.e. a homeomorphism to a model topological surface, with coordinates obtained from a pair of pants decomposition by lifting the angle $\theta({\rm mod\ }l)$ to take values in $\br$.  The lengths $l_i\in\br^+$ and angles $\theta_i\in\br$ give global coordinates called {\em Fenchel-Nielsen coordinates} on the Teichm\"uller space of marked oriented hyperbolic surfaces.
The orientation is required to make sense of the angle at each two-sided geodesic.  One can choose at each two-sided curve the way in which to define $\theta$ (for example, from one side of a geodesic in the decomposition turn right at a distinguished point and head towards a point distinguished from the other side) to leave an ambiguity $\theta\mapsto\theta+l/2({\rm mod\ }l)$ compared to (\ref{eq:theta}). Furthermore, the change from $\theta$ to $\theta+l/2$ gives a different hyperbolic surface so we are left with $\theta({\rm mod\ }l)$.  A flip $\theta\mapsto l-\theta$ at every two-sided geodesic gives the same hyperbolic surface with the opposite orientation so we have a double cover of the moduli space of unoriented surfaces.  

Such a double cover of the moduli space requires more care when the surface is non-orientable.  Teichm\"uller spaces of non-orientable hyperbolic surfaces have been studied by various authors.  In \cite{JosMin, SepTei} Fenchel-Nielsen coordinates are defined on marked non-orientable hyperbolic surfaces.  Jost \cite{JosMin} considers an enlarged Teichm\"uller space that combines orientable and non-orientable surfaces which naturally arises when compactifying the spaces.  The fact that one must deal with some type of an orientation on the surfaces (see below) is not explicitly mentioned in these papers although it is implicit in \cite{JosMin}.   

We define the ``oriented" moduli space of surfaces when the boundary is non-empty.   
\begin{defn}
For $n>0$ and $n+k-2>0$, denote by $\widehat{\modm}_{k,n,\pm}({\bf L})$ the moduli space of hyperbolic surfaces with $n$ geodesic boundary components $\beta_1,...,\beta_n$ of lengths ${\bf L}=(L_1,...,L_n)$ and fixed topology, together with an orientation of $\beta_1$.
\end{defn}
In the orientable case $\widehat{\modm}_{k,n,+}({\bf L})=\modm_{k/2,n}({\bf L})$ agrees with the usual definition.   In the non-orientable case, on this oriented moduli space the volume form is well-defined without the $\pm 1$ ambiguity present in Theorem~\ref{th:vol}.  This can be seen by choosing a curve $\gamma$ with connected orientable complement unique up to the action of the mapping class group.   The orientation of $\beta_1$ induces an orientation on the complement of $\gamma$.  Upstairs, in Teichm\"uller space, use coordinates coming from a decomposition along a maximal set of disjoint embedded closed curves that includes a curve $\gamma$ with connected orientable complement.  The orientation on the complement of $\gamma$ gives rise to a well-defined angle at each of the closed geodesics in the decomposition and hence a well-defined volume form.  We describe this explicitly below in the simplest case of a Klein bottle minus a disk in order to demonstrate an integration formula.  More generally, for any decomposition of a surface $\Sigma$ along a maximal set of disjoint embedded closed curves, a well-defined orientation, which depends on an order of the one-sided geodesics in the maximal set, is induced from an orientation on $\beta_1$ as follows.  Cut along the one-sided geodesics in the maximal set to get an orientable surface $\Sigma'$.  The orientation on $\beta_1$ induces an orientation on $\Sigma'$.  To determine the order of the boundary one-sided geodesics, join pairs of one-sided geodesics by a collection of disjoint embedded arcs in $\Sigma'$.  (There may be a one-sided geodesic left over.)  The two one-sided geodesics can be replaced by a two-sided closed geodesic in $\Sigma$ which restricts to the isotopy class of the arc in $\Sigma'$.  The choice of angle along the two-sided closed geodesic is well-defined, induced from the orientation on $\beta_1$.  We then use (\ref{eq:symp2}) which chooses an order of the two one-sided geodesics that makes the local volume form agree.

The oriented moduli space can also be described by an index two subgroup of the mapping class group.  Take those elements of the mapping class group that preserve the chosen orientation on the specified boundary component.  To prove that this does indeed give an index two subgroup (and not the whole group) we must show that there exists an element of the mapping class group that flips the orientation.  
\begin{lemma}
On any surface with boundary there exists a homeomorphism that flips a chosen orientation on a boundary component.
\end{lemma}
\begin{proof}
The simplest cases are:  a M\"obius strip, a M\"obius strip minus a disk, and an orientable surface.  A pair of pants possesses an orientation reversing reflection and this extends to any orientable surface cut along pairs of pants.  The disk and annulus also possess reflections.  For the M\"obius strip minus a disk the mapping class group is $\bz_2$ and it is easy to see explicitly that the non-trivial element flips the orientation of both boundary components.  Furthermore, this homeomorphism extends across the disk to the M\"obius strip.

In general, any non-orientable surface $\Sigma$ with boundary contains a M\"obius strip minus a disk $M$ with one boundary component in common with a boundary component of $\Sigma$.  A homeomorphism of $M$ that flips the orientation of the boundary extends over all of $\Sigma$ by choosing a homeomorphism of $\Sigma-M$ that flips the orientation of its boundary component shared with $M$.  This homeomorphism of $\Sigma-M$ exists by induction on the simpler surface $\Sigma-M$.  Thus the lemma is proven.
\end{proof}
The kernel of the action of the mapping class group on the orientation on a specified boundary component gives the following corollary.
\begin{cor}
There exists an index two subgroup of the mapping class group of a non-orientable surface with boundary that fixes a chosen orientation on a boundary component thus defining a double cover of the moduli space.
\end{cor}
It is not clear that this index two subgroup exists for a closed non-orientable surface.  Equivalently,  it is not clear that there exists a subgroup of the mapping class group that preserves the volume form which is well-defined on the Teichm\"uller space of a closed non-orientable surface.  

In place of the decomposition of a hyperbolic surface along a maximal set of disjoint embedded closed geodesic it is also convenient to use Penner's coordinates \cite{PenDec} of Teichm\"uller space given by lengths of a maximal set of disjoint geodesic arcs perpendicular to the boundary.  These coordinates consist of lengths and no angles.   Lengths are a little more natural on a non-oriented surface.  In this case the $\pm1$ ambiguity of the volume form on the moduli space comes from a choice of ordering for the coordinates just as for lengths of one-sided geodesics.  In \cite{MonTri} the Poisson structure on Teichm\"uller space is given in these coordinates.  This Poisson structure makes sense on the moduli spaces of non-orientable surfaces, and is quite natural particularly since the dimensions of the moduli spaces can be odd.  The Poisson structure requires a choice of orientation on each piece of the marked hyperbolic surface cut along the perpendicular geodesics.  This can be achieved in the same way as for the volume form.

\subsection{Integration}
Mirzakhani showed how to use McShane's identity, and her generalised versions of the identity, to integrate functions on the moduli space, and in particular calculate the volume of the moduli space \cite{MirGro,MirSim,MirWei}. 

We give a brief description of her idea, before we apply it to a non-orientable example.  Mirzakhani shows how to integrate functions of the form 
\[ F=\sum_{\gamma=h\cdot\gamma_0} f(l_{\gamma})\] 
where $f$ is an arbitrary function and the sum is over the orbit of a geodesic under the mapping class group.  When $f$ decays fast enough the sum is well-defined on the moduli space.  More generally,
one can consider an arbitrary (decaying) function on collections of geodesics and sum over orbits of the mapping class group acting on the collection.   Mirzakhani unfolds the integral of $F$ to a moduli space $\widehat{\modm}_{g,n}({\bf L})$ of pairs $(\Sigma,\gamma)$ consisting of a hyperbolic surface $\Sigma$ and a geodesic $\gamma\subset\Sigma$.  
\[\begin{array}{c}
\modt_{g,n}({\bf L})\\\downarrow\\\widehat{\modm}_{g,n}({\bf L})\\\downarrow\\\modm_{g,n}({\bf L})
\end{array}\]
The unfolded integral 
\[\int_{\modm_{g,n}({\bf L})}F\cdot d{\rm vol}=\int_{\widehat{\modm}_{g,n}({\bf L})}f(l_{\gamma})\cdot d{\rm vol}\]
can be expressed in terms of an integral over the simpler moduli space obtained by cutting $\Sigma$ along the geodesic $\gamma$.  

McShane's identity is exactly of the right form for Mirzakhani's scheme since it expresses the constant function $F=L_1$ as a sum of functions of lengths over orbits of the mapping class group.  In this case,
\[ L_1V_{g,n}({\bf L})=\int_{\modm_{g,n}({\bf L})}F\cdot d{\rm vol}=\int_{\widehat{\modm}_{g,n}({\bf L})}f(l_{\gamma_1},l_{\gamma_2})\cdot d{\rm vol}\]
expresses the volume $V_{g,n}({\bf L})=$ volume$(\modm_{g,n})$ recursively in terms of simpler volumes of moduli spaces obtained by removing a pair of pants from $\Sigma$ in topologically different ways.  The simplest unfolding expresses the volume of $\modm_{1,1}$ as an integral over $\br^+\times[-1/2,1/2]$.  In the example of the Klein bottle below, we see a similar unfolding to an integral over $\br^+\times\br^+$.

The volume of $\modm_{k,n,-}({\bf L})$ is infinite so we cannot integrate constant functions over the moduli space but we can integrate integrable functions.  

As described above, we can take a double cover of $\modm_{2,1,-}$, the moduli space of punctured Klein bottles (we have dropped the length $L$ since it is zero), on which the volume form is well-defined.  The double cover is obtained by orienting the boundary which induces an orientation on a complement of the two-sided geodesic $\gamma$.  
The volume form is 
\[ dl_{\gamma}\wedge d\theta_{\gamma}=\coth\frac{1}{2}l_{\gamma_i}\coth\frac{1}{2}l_{\gamma_{i+1}}dl_{\gamma_i}\wedge dl_{\gamma_{i+1}}\]
since the orientation has allowed us to choose the sign of $\theta_{\gamma}$, and also to distinguish $\gamma_{i-1}$ from $\gamma_{i+1}$ with respect to $\gamma_i$.

The length $l_{\gamma}$ is well defined on the moduli space since $\gamma$ is unique.  The identity from Theorem~\ref{th:kident}
\[ 
\sum_{i=-\infty}^{\infty}\frac{1}{1+\sinh^2\frac{1}{2}l_{\gamma_i}+\sinh^2\frac{1}{2}l_{\gamma_{i+1}}}=\tanh\frac{1}{2}l_{\gamma}
\]
enables us to unfold the integral of $\tanh\frac{1}{2}l_{\gamma}$ over the moduli space, but this integral is infinite.  Instead, integrate
\[ F=\frac{\tanh\frac{1}{2}l_{\gamma}}{\cosh^n\frac{1}{2}l_{\gamma}}=
\sum_{i=-\infty}^{\infty}\frac{(2\sinh\frac{1}{2}l_{\gamma_i}\sinh\frac{1}{2}l_{\gamma_{i+1}})^n}{(1+\sinh^2\frac{1}{2}l_{\gamma_i}+\sinh^2\frac{1}{2}l_{\gamma_{i+1}})^{n+1}}.\]
This identity uses (\ref{eq:quad}) and the identity for $n=0$.

We will calculate only the simplest case $n=1$
\begin{eqnarray*}
\int_{\modm_{2,1,-}}F\cdot d{\rm vol}
&=&\int_{\modm_{2,1,-}}\sum_{i=-\infty}^{\infty}\frac{2\sinh\frac{1}{2}l_{\gamma_i}\sinh\frac{1}{2}l_{\gamma_{i+1}}}{(1+\sinh^2\frac{1}{2}l_{\gamma_i}+\sinh^2\frac{1}{2}l_{\gamma_{i+1}})^2}d{\rm vol}\\
&=&\int_{\widehat{\modm}_{2,1,-}}\frac{2\sinh\frac{1}{2}l_{\gamma_i}\sinh\frac{1}{2}l_{\gamma_{i+1}}\coth\frac{1}{2}l_{\gamma_i}\coth\frac{1}{2}l_{\gamma_{i+1}}}{(1+\sinh^2\frac{1}{2}l_{\gamma_i}+\sinh^2\frac{1}{2}l_{\gamma_{i+1}})^2}
dl_{\gamma_i}dl_{\gamma_{i+1}}\\
&=&\int _0^{\infty}\int _0^{\infty}\frac{2\cosh{\frac{1}{2}x}\cosh{\frac{1}{2}y}}{(1+\sinh^2{\frac{1}{2}x}+\sinh^2{\frac{1}{2}y})^2}dxdy.\\
&=&2\pi\end{eqnarray*}
In this case, the integral can be performed without unfolding, since
\[ \int _0^{\infty}\int_0^{l_{\gamma}}\frac{\tanh\frac{1}{2}l_{\gamma}}{\cosh\frac{1}{2}l_{\gamma}}d\theta_{\gamma} dl_{\gamma}=2\pi.\]

Mirzakhani \cite{MirWei} used the volumes of the moduli spaces to study cohomology classes on the moduli space.  She proved recursion relations between volumes of the moduli spaces which led to recursion relations between cohomology classes, reproving the Witten-Kontsevich Virasoro relations on the moduli space.  It is not clear what the analogous setup should be on the moduli spaces of non-orientable surfaces since the volumes are infinite.  Nevertheless, one would hope to be able to integrate over the moduli spaces of non-orientable surfaces and get relations between cohomology classes on different moduli spaces.  
The degree $4i$ tautological cohomology classes appearing in Wahl's work \cite{WahHom} on the cohomology of moduli spaces of non-orientable surfaces are the natural candidates to satisfy recursion relations.

\end{document}